\documentclass[11pt]{article}

\newcommand{\junk}[1]{}

\usepackage[left=1in,right=1in,top=1in,bottom=1in]{geometry}
\usepackage{amsmath,amssymb,fancyhdr}
\usepackage{enumerate}
\usepackage{mdframed}
\usepackage{lineno}
\usepackage{hyperref}

\usepackage{epsfig}
\usepackage{amsfonts}
\usepackage{amssymb}
\usepackage{amstext}
\usepackage{xspace}
\usepackage{theorem}
\usepackage{graphicx}

\usepackage{algorithmicx}
\usepackage{algpseudocode}

\theoremstyle{definition}
\newtheorem{theorem}{Theorem}

\newtheorem{proposition}[theorem]{Proposition}

\newcommand{\Z}{\mathbb{Z}}

\newenvironment{proofof}[1]{{\bf Proof of #1:  }}{\hfill\rule{2mm}{2mm}}

\title {Linear $d$-polychromatic $Q_{d-1}$-colorings of the Hypercube}

\author{
Eugene Han\\
\small Department of Mathematical Sciences\\[-0.8ex]
\small Carnegie Mellon University\\[-0.8ex]
\small Pittsburgh, PA, USA\\
\small\tt eugeneh@andrew.cmu.edu
\and
David Offner\\
\small Department of Mathematics and Computer Science\\[-0.8ex]
\small Westminster College\\[-0.8ex] 
\small New Wilmington, PA, USA\\
\small\tt offnerde@westminster.edu\\
}

\begin{document}

\maketitle

\begin{abstract}
Let $n \ge d \ge \ell \ge 1$ be integers, and denote the $n$-dimensional hypercube by $Q_n$. A coloring of the $\ell$-dimensional subcubes $Q_\ell$ in $Q_n$ is called a $Q_\ell$-coloring.  Such a coloring is $d$-polychromatic if every $Q_d$ in the $Q_n$ contains a $Q_\ell$ of every color.  In this paper we consider a specific class of $Q_\ell$-colorings that are called linear. Given $\ell$ and $d$, let $p_{lin}^\ell(d)$ be the largest number of colors such that there is a $d$-polychromatic linear $Q_\ell$-coloring of $Q_n$ for all $n \ge d$.  We prove that for all $d \ge 3$, $p_{lin}^{d-1}(d) = 2$. In addition, using a computer search, we determine $p_{lin}^\ell(d)$ for some specific values of $\ell$ and $d$, in some cases improving on previously known lower bounds.
\end{abstract}

\section{Introduction}

For $n \in \Z$, $n \ge 1$,  the \emph{$n$-dimensional hypercube}, denoted
by $Q_n$, is the graph with $V(Q_n) = \{0,1\}^n$, and edges
between vertices which differ in exactly one coordinate. For $\ell \le n$, $Q_\ell$-coloring of $Q_n$ is a coloring of the $\ell$-dimensional subcubes $Q_\ell$ in $Q_n$.  A \textit{$Q_\ell$-coloring} of $Q_n$ is called \textit{$d$-polychromatic} if every $Q_d$ in $Q_n$ contains a $Q_\ell$ of every color. For $1 \le \ell \le d$, let $p^\ell(d)$ be the maximum number of colors such that for all $n \ge d$ there is a $d$-polychromatic $Q_\ell$-coloring of $Q_n$.

Motivated by Tur\'an type problems on the hypercube, Alon, Krech, and Szab\'o \cite{AKS07} introduced the notion of polychromatic colorings on the hypercube and proved bounds on $p^1(d)$.
\begin{theorem}[Alon, Krech, and Szab\'o~\cite{AKS07}]\label{AKS}
For all $d \ge 1$,
\[ \binom{d+1}{2} \ge p^1(d) \ge \left\lfloor \frac{(d+1)^2}{4} \right\rfloor. \]
\end{theorem}

The lower bound is given by what Chen~\cite{Chen16} called \textit{basic colorings}, which we define in Section~\ref{back}.  In \cite{Off08}, it was shown that for edge colorings of the hypercube, a basic coloring is optimal.
\begin{theorem}[Offner~\cite{Off08}]\label{pd}
For all $d \ge 1$,
  \[ p^1(d) = \left\lfloor \frac{(d+1)^2}{4} \right\rfloor. \]
\end{theorem}

Later, \"Ozkahya and Stanton generalized the bounds of Theorem~\ref{AKS} to $p^\ell(d)$ for $\ell >1$.
\begin{theorem}[\"Ozkahya and Stanton~\cite{OS11}]\label{OS}
For all $d \ge \ell \ge 1$, let $0 <r \le \ell + 1$ be such that $r \equiv d+1 \pmod{\ell+1}$. Then
\[ \binom{d+1}{\ell + 1} \ge p^\ell(d) \ge \left\lceil \frac{d+1}{\ell + 1} \right\rceil^{r}\left\lfloor \frac{d+1}{\ell + 1} \right\rfloor^{\ell + 1-r}. \]
\end{theorem}

As with Theorem~\ref{AKS}, the lower bound in Theorem~\ref{OS} is achieved by basic colorings.
Define $p_{bas}^\ell(d)$  to be the maximum number of colors such that any hypercube has a $d$-polychromatic basic $Q_\ell$-coloring. Then, as described in Section~\ref{back},

\begin{equation}\label{pldbas}
p_{bas}^\ell(d) = \left\lceil \frac{d+1}{\ell + 1} \right\rceil^{r}\left\lfloor \frac{d+1}{\ell + 1} \right\rfloor^{\ell + 1-r}. 
\end{equation}

In contrast to the case of $Q_1$-colorings, where for all $d\ge 1$, $p^1(d) = p^1_{bas}(d)$, it is not true in general that $p^\ell(d) = p^\ell_{bas}(d)$.  This was first shown in $\cite{GLMOTY16}$.

\begin{theorem}[\cite{GLMOTY16}]\label{GLMOTY}
$p^2(3) = 3$ and $p^2(4) \ge 5$.%NOTE: cite authors?

\end{theorem}
Note that $p_{bas}^2(3) = 2<p^2(3)$ and $p_{bas}^2(4) = 4 < p^2(4)$.  Subsequently, Chen~\cite{Chen16} improved the lower bound on $p^2(d)$ for all $d \ge 4$.

\begin{theorem}[Chen~\cite{Chen16}]\label{Chen}
For all $d \ge 4$, 
\[ p^2(d) \ge \begin{cases}
(k^2+1)(k+1) & d=3k;\\
(k^2+k+1)(k+1) & d=3k+1;\\
(k^2+k+1)(k+2) & d=3k+2.\\
\end{cases}
 \]
\end{theorem}

A corollary of Theorem~\ref{GLMOTY} (See \cite{GLMOTY16}, Corollary 24) is that for all $d \ge 3$, $p^{d-1}(d) > p_{bas}^{d-1}(d)$, and a corollary of Theorem~\ref{Chen} (See \cite{Chen16}, Corollary 1.3) is that for all $\ell >1$, $d > \ell+1$, $p^\ell(d) > p_{bas}^\ell(d)$. Thus for all $d > \ell > 1$, $p^\ell(d) > p_{bas}^\ell(d)$. 

Chen called the colorings used to establish the bounds in Theorem~\ref{Chen} \textit{linear colorings}, a term defined in Section~\ref{back}. Define $p_{lin}^\ell(d)$ to be the maximum number of colors such that for all $n \ge d$, $Q_n$ has a $d$-polychromatic linear $Q_\ell$-coloring.  All basic colorings are linear, so for all $d\ge \ell \ge 1$, 

\[p_{bas}^\ell(d) \le p_{lin}^\ell(d) \le p^\ell(d).\]

For $d \ge 4$, the best upper bound on $p^2(d)$ is from Theorem~\ref{OS}, $p^2(d) \le \binom{d+1}{3}$.
 In~\cite{Chen16}, Chen used a geometric argument to prove a nontrivial upper bound on $p_{lin}^2(d)$.
\begin{theorem}[Chen~\cite{Chen16}]\label{ch16}
For all $d$ sufficiently large,
\[p_{lin}^2(d) \le \frac{26}{27}\binom{d+1}{3}.\] 
\end{theorem}

Since the colorings used in proving Theorem~\ref{Chen} are linear, we know that for all $\ell > 1$, $d > \ell +1$, $p_{bas}^{\ell}(d) < p_{lin}^{\ell}(d)$.  It is natural to ask whether for $d \ge 3$, $p_{bas}^{d-1}(d)<p_{lin}^{d-1}(d)$.  The main result of this paper shows that the answer to this question is no.

\begin{theorem}\label{linnotopt}
For all $d \ge 3$, $p_{lin}^{d-1}(d) = 2$.
\end{theorem}

Section~\ref{pfmain} is devoted to the proof of Theorem~\ref{linnotopt}.

To determine the values of $p_{lin}^\ell(d)$ for some particular values of $\ell$ and $d$, we wrote a program that, given $\ell$, $d$, and $M$, tests whether a $d$-polychromatic linear $Q_\ell$-coloring with $M$ colors exists for all $n$.  The program is described in Section~\ref{code} and the python code is available at \texttt{arXiv.org}. Theorem~\ref{Chen} and its corollary imply the following lower bounds:
\[\begin{matrix}
p_{lin}^2(4) \ge 6, & p_{lin}^2(5) \ge 9, & p_{lin}^3(5) \ge 6, & p_{lin}^2(6) \ge 15, & p_{lin}^3(6) \ge 9, &p_{lin}^4(6) \ge 6.
\end{matrix}\]
In the first three of these cases our search did not turn up any polychromatic linear colorings with more colors, proving that those linear colorings are optimal.  In the other cases we found linear colorings which establish larger values for $p^\ell_{lin}(d)$. These results are summarized in the next proposition.

\begin{proposition}\label{small}
$p_{lin}^2(4) = 6$, $p_{lin}^2(5) = 9$, $p_{lin}^3(5) = 6$, $p_{lin}^2(6) = 16$, $p_{lin}^3(6) = 12$, and $p_{lin}^4(6) = 9$.
\end{proposition}
We prove Proposition~\ref{small} in Section~\ref{dere}.

For $d \ge 7$, the problem of determining $p_{lin}^\ell(d)$ for $2 \le \ell \le d-2$ remains open.  The colorings used in the last three cases of Proposition~\ref{small} are the only known polychromatic linear colorings with more colors than the linear colorings used in the proof of Theorem~\ref{Chen}.  It would be interesting to further improve the lower bounds in Theorem~\ref{Chen}, or to improve the upper bound on linear colorings in Theorem~\ref{ch16}.

Since $p_{lin}^2(4) = p_{lin}^3(5) = 6$ and $p_{lin}^4(6) = 9$, unlike the case of $p_{lin}^{d-1}(d)$, it is not true that for all $d \ge 4$, $p_{lin}^{d-2}(d)$ is equal to a constant.  It would be interesting to determine these values.

For general polychromatic colorings, the only known value of $p^\ell(d)$  where $2 \le \ell < d$ is $p^2(3) = 3$. Combining Theorem~\ref{linnotopt} with Theorems~\ref{GLMOTY} and \ref{OS}, we conclude that for all $d \ge 3$, 
\[2=p_{bas}^{d-1}(d)=p_{lin}^{d-1}(d) < 3 \le p^{d-1}(d) \le d+1.\]
These are the only examples where it is known that $p_{lin}^\ell(d) < p^\ell(d)$. Additionally, Theorems~\ref{OS} and \ref{Chen} imply that for all $d \ge 4$,
\[ \binom{d+1}{3} \ge p^2 (d) \ge \begin{cases} 
  (k^2+1)(k+1) & d=3k \\ 
  (k^2+k+1)(k+1) & d=3k +1\\ 
    (k^2+k+1)(k+2) & d=3k +2. 
  \end{cases} \]
It would be interesting to prove whether there are other cases where $p_{lin}^\ell(d) < p^\ell(d)$ or improve any of these bounds on $p^\ell(d)$.

\section{Basic and Linear Colorings}\label{back}

We  follow standard notation and denote $d$-dimensional subcubes of $Q_n$ by an $n$-bit binary string with $d$ entries replaced by stars, where each vertex of the $Q_d$ is obtained by replacing the stars with 0's and 1's.  For example, $[01{*}10{*}011]$ is a copy of $Q_2$ in $Q_9$ with vertices $[010100011]$, $[011100011]$, $[010101011]$, and  $[011101011]$.

Given a copy of $Q_d$ in $Q_n$, we define its \textit{counting vector} $(v_0,\ldots, v_d) \in \Z^{d+1}_{\ge 0}$, to be the vector where $v_0$ is the number of 1's to the left of the first star, $v_d$ is the number of 1's to the right of the $d$th star, and for $1 \le i \le d-1$, $v_i$ is the number of 1's between the $i$th and $(i+1)$st stars. For example the $Q_3$'s $[01{*}1011{*}00{*}1010]$ and $[1{*}011100{*}{*}00110]$ both have the counting vector $(1,3,0,2)$. A $Q_\ell$-coloring $\chi$ is called \textit{simple} if all $Q_\ell$'s with the same counting vector are assigned the same color.  An application of Ramsey's theorem (See e.g. Lemma 3 from~\cite{GLMOTY16}) 
implies that when studying polychromatic colorings on the hypercube, we need only consider simple colorings, so all $Q_\ell$-colorings in this paper will have the form
\[\chi: \Z^{\ell+1}_{\ge 0} \to S,\]
where $S$ is some set of colors, and $\chi$ colors any $Q_\ell$ with the color of its counting vector.

 A $Q_\ell$-coloring $\chi$ is called \textit{linear} if the set of colors is a finite abelian group $Z$, and the coloring is induced by an additive map 
\[\chi: \Z^{\ell+1}_{\ge 0} \to Z\]

A linear $Q_\ell$-coloring $\chi$ is called $(\ell,d)$-\textit{basic} (or just \textit{basic} when the values of $\ell$ and $d$ are clear from context) when the finite abelian group $Z$ is of the form $Z=\bigoplus^{\ell}_{i=0} \Z/m_i$, where $m_0 + \cdots + m_\ell =d+1$, and 
\[\chi(v_0, \ldots, v_{\ell})=(v_0\pmod{m_0}, \ldots, v_{\ell}\pmod{m_{\ell}}).\]
Note that a $(\ell,d)$-basic coloring uses $m_0m_1\cdots m_\ell$ colors, and in \cite{OS11} it is shown that a $(\ell,d)$-basic coloring is always $d$-polychromatic. The value of $p_{bas}^\ell(d)$ from Equation~\ref{pldbas} is obtained by choosing the values of $m_i$ so that for all $0 \le i<j \le \ell$, $|m_i-m_j|\le 1$.

\section{Linear $d$-polychromatic $Q_{d-1}$-colorings have at most 2 colors}\label{pfmain}

\begin{proofof}{Theorem~\ref{linnotopt}} Fix $d \ge 3$ and consider a linear $Q_{d-1}$-coloring $\chi:\Z^{d}_{\ge 0} \to Z$ of a hypercube $Q_n$, where $n$ is very large, $Z$ is a finite abelian group, and $|Z| \ge 3$.  Our goal is to show that $\chi$ is not $d$-polychromatic. For $0 \le i \le d-1$, let $f_i \in \Z^{d}_{\ge 0}$ be the vector whose $(i+1)$th entry is a 1, and all other entries are 0, and let $h_i \in Z$ be such that 
\[h_i = \chi(f_i).\]
Then for any $(v_0, v_1, \ldots, v_{d-1}) \in \Z^{d}_{\geq 0}$, 
\[ \chi(v_0, v_1, \ldots, v_{d-1}) = \sum\limits^{d-1}_{i=0}v_{i} h_{i}.\]

In this proof, let $(x_0, x_1, \ldots, x_d)$ be the counting vector for a $Q_d$. To show $\chi$ is not $d$-polychromatic, we examine the counting vectors of the $2d$ $Q_{d-1}$'s contained in this $Q_d$. For $0 \le i \le d-1$, let $v_{i0}$ and $v_{i1}$ denote the counting vectors of the two $Q_{d-1}$'s in such a $Q_d$ that use all but the $(i+1)$th star, i.e.
\[ v_{i0} = (x_0, \ldots, x_{i-1}, x_{i}+ x_{i+1}, x_{i+2}, \ldots, x_d), v_{i1} =(x_0, \ldots, x_{i-1}, x_{i}+ x_{i+1}+1, x_{i+2}, \ldots, x_d).\]
For example, if $d=4$, $\ell=3$, and we consider the $Q_4$ represented by $[0100{*}11{*}0{*}11{*}01]$, then $(x_0, x_1, x_2, x_3, x_4)= (1,2,0,2,1)$, and $v_{00} = (3,0,2,1)$, $v_{01} = (4,0,2,1)$, $v_{10} = (1,2,2,1)$, $v_{11} = (1,3,2,1)$, $v_{20} = (1,2,2,1)$, $v_{21} = (1,2,3,1)$, $v_{30} = (1,2,0,3)$, and $v_{31} = (1,2,0,4)$.

For $0 \le i \le d-1$, let $c_i$ denote the set of colors assigned to those two $Q_{d-1}$'s that use all but the $(i+1)$th star, i.e. \[c_i = \{\chi(v_{i0}), \chi(v_{i1}),\}\]
and let \[C_i = \bigcup_{j=0}^i c_j.\]
Note that $C_{d-1}$ represents the set of colors assigned to all $Q_{d-1}$'s in such a $Q_d$. 

Let 
\[X = \chi(v_{00}) = x_0h_0 + x_1h_0 + x_2h_1 + \cdots + x_{d}h_{d-1}.\] 
Note that for $i \ge 0$, \[v_{(i+1)0} - v_{i0}  = x_{i+1}(f_{i+1}-f_i) =  (0,\ldots,0, -x_{i+1},x_{i+1},0,\ldots,0).\]
Thus if we let $H_0=0$, and for all $1\le i \le d-1$, define
\[H_i = \sum_{j=1}^i x_{j}(h_j - h_{j-1}),\]
then for  $0 \le i \le d-1$, 
 \[c_i  = \{X + H_i, X + H_i + h_i\}.\]
Further, note that $H_i$ is fixed by the choices of $x_0, \ldots, x_{i}$, and for $i \ge 0$, $H_{i+1} = H_i + x_{i+1}(h_{i+1} - h_{i})$.

%We prove that $\chi$ is not $d$-polychromatic.  To do this, we show that there are choices for $x_0, x_1, \ldots, x_{d+1} \in \Z_{\geq 0}$ such that $|C_{d+1}| = |\bigcup c(F_i) | < |G|$. 

If $h_i=0$ for all $i$, then $\chi$ colors every $Q_{d-1}$ with the color 0, and is clearly not $d$-polychromatic.  Thus we assume $h_i \neq 0$ for some $i$. Let $0 \le i_1<i_2 < \cdots < i_k \le d-1$ be the indices such that for $1 \le j \le k$, $h_{i_j} \neq 0$.

In this proof, we show that there are choices for $x_0, x_1, \ldots, x_{d} \in \Z_{\ge 0}$ such that if $(x_0, x_1, \ldots, x_d)$ is the counting vector for a $Q_d$, it will not contain a $Q_{d-1}$ of each color. Let $g \in Z$ such that $g \notin \{0, h_{i_1}\}$. We prove the following statement by induction:  There are choices for $x_0, \ldots, x_{{i_k}} \in \Z_{\ge 0}$ such that for any choice of $x_{(i_k)+1}, \ldots, x_{d}$, $X+g \notin C_{i_k}$.

\textbf{Base case $(j=1)$:} Set $x_0=x_1 = \cdots =x_{i_1} = 0$, so  $H_0 = H_1 = \cdots  = H_{i_1} = 0$.  Then $c_{i_1} = \{X, X+h_{i_1}\}$, and if $i_1>0$, since $h_i = 0$ for $i < i_1$, $C_0 = C_1 = \cdots = C_{i_1-1} = \{X\}$. Thus $C_{i_1} = \{X, X+h_{i_1}\}$. Since $g \notin \{0, h_{i_1}\}$, $X+g \notin C_{i_1}$.

\textbf{Inductive step:} Fix $j$ such that $1 \le j \le k-1$ and assume there are choices for $x_0, \ldots, x_{i_{j}}$ such that for any choices of $x_{(i_j)+1}, \ldots, x_{d}$, $X+g \notin C_{i_{j}}$.
We want to show there are choices for $x_{(i_j)+1}, \ldots, x_{i_{j+1}}$ such that for any choice of $x_{(i_{j+1})+1}, \ldots, x_{d}$, $X+g \notin C_{i_{j+1}}$. 

First note that 
\[c_{i_j}= \{X + H_{i_j}, X + H_{i_j} + h_{i_j}\},\] 
so by the induction hypothesis, $g \neq H_{i_j}$ and $g \neq H_{i_j} + h_{i_j}$.
We treat two cases, depending on whether $h_{(i_j)+1}$ is 0 or not.

\textbf{Case 1:} Suppose $h_{(i_j)+1} \neq 0$, i.e. $i_{j+1} = (i_j)+1$. We show that setting $x_{(i_j)+1}=0$, $|Z|-1$, or 1 ensures that $X+g \notin C_{i_{j+1}}$.

\textbf{Case 1a:} Suppose $g \neq H_{i_j} +h_{(i_j)+1}$. In this case, let $x_{(i_j)+1}=0$. Then $H_{(i_j)+1} = H_{i_j}$, so 
\[c_{i_{j+1}} = c_{(i_j)+1} =  \{X +H_{(i_j)+1} , X +H_{(i_j)+1} +h_{(i_j)+1}\} = \{X +H_{i_j} , X +H_{i_j} +h_{(i_j)+1}\},\] 
and this choice works since $g \neq H_{i_j}+h_{(i_j)+1}$.

\textbf{Case 1b:} Suppose $g \neq H_{i_j}+ h_{i_j}-h_{(i_j)+1}.$ In this case, let $x_{(i_j)+1}= |Z|-1$.  Then $H_{(i_j)+1} = H_{i_j} -(h_{(i_j)+1} - h_{i_j})$, so
\[\begin{split}
c_{i_{j+1}} &= c_{(i_j)+1} \\
&=   \{X+H_{i_j} -(h_{(i_j)+1} - h_{i_j}), X+H_{i_j} -(h_{(i_j)+1} - h_{i_j})+h_{(i_j)+1}\}\\
&=\{X+H_{i_j} -(h_{(i_j)+1} - h_{i_j}), X+H_{i_j}+ h_{i_j})\} ,
\end{split}\]
and this choice works since  $g\neq H_{i_j} -(h_{(i_j)+1} - h_{i_j})$.

\textbf{Case 1c:} The only remaining case occurs if we are not in Case 1a, so
\begin{equation}\label{gmh1}
g-H_{i_j}=h_{(i_j)+1},
\end{equation} and we are not in Case 1b, so
\begin{equation}\label{gmh2}
g-H_{i_j} = h_{i_j}-h_{(i_j)+1}.
\end{equation}
Combining Equations~\ref{gmh1} and \ref{gmh2}, we get 
\begin{equation}\label{2hj+1hj}
2h_{(i_j)+1} = h_{i_j}.
\end{equation}

In this case, set $x_{(i_j)+1}=1$.  Then
\[\begin{split} c_{i_{j+1}} &= c_{(i_j)+1}    \\
&=\{X+H_{i_j} +(h_{(i_j)+1} - h_{i_j}), X+H_{i_j} +(h_{(i_j)+1} - h_{i_j}) +h_{(i_j)+1}\}\\
&=\{X+H_{i_j} +(h_{(i_j)+1} - 2h_{(i_j)+1}), X+H_{i_j} +(h_{(i_j)+1} - 2h_{(i_j)+1}) +h_{(i_j)+1}\}\\
&=\{X+H_{i_j} -h_{(i_j)+1}, X+H_{i_j}\}.
\end{split}\]
%where the last equality follows from Equation~\ref{2hj+1hj}.

This choice works unless $g=H_{i_j} -h_{(i_j)+1}$, which implies
\begin{equation}\label{gmh3}
g-H_{i_j} =-h_{(i_j)+1}.
\end{equation}

Equations~\ref{gmh2} and \ref{gmh3} together imply that $h_{i_j}=0$, a contradiction. %Thus one of the three values $x_{i_j}=0, |Z|-1, 1$ works.

\textbf{Case 2:} Suppose $h_{(i_j)+1} = 0$, i.e. $i_{j+1} > (i_j)+1$. We show that setting $x_{(i_j)+1} = 0$ or $|Z|-1$ along with $x_{(i_j)+2} = x_{(i_j)+3} = \cdots = x_{i_{j+1}}=0$ ensures that $X+g \notin C_{i_{j+1}}$.

Since $h_{(i_j)+1} = 0$,
\[\begin{split}
c_{(i_j)+1} &=\{X + H_{(i_j)+1}, X + H_{(i_j)+1} + h_{(i_j)+1}\}\\
&= \{X + H_{i_j} + x_{(i_j)+1}(h_{(i_j)+1}-h_{i_j}), X + H_{i_j} + x_{(i_j)+1}(h_{(i_j)+1}-h_{i_j}) + h_{(i_j)+1}\}\\
&=\{X + H_{i_j} -x_{(i_j)+1}h_{i_j}\}.
\end{split}\]

Thus since $g \neq H_{i_j}$ and $g \neq H_{i_j} + h_{i_j}$, a choice of $x_{(i_j)+1} = 0$ or $x_{(i_j)+1} = |Z|-1$ will ensure $X + g \notin C_{(i_j)+1}$.  

For $(i_j)+2 \le i \le (i_{j+1})-1$, since $h_i = h_{i-1}= 0$ we know $H_i=H_{i-1}$ no matter the choice of $x_{i}$.  Thus any values of $x_i$ for $(i_j)+2 \le i \le (i_{j+1})-1$ ensures that $X + g \notin C_{i}$ for $(i_j)+2 \le i \le (i_{j+1})-1$.

Finally, since $ H_{(i_{j+1})-1} = H_{i_j} -x_{(i_j)+1}h_{i_j}$, by setting $x_{i_{j+1}}=0$ we get
\[H_{i_{j+1}} = H_{(i_{j+1})-1}+ x_{i_{j+1}}(h_{i_{j+1}} - h_{(i_{j+1})-1}) = H_{i_j} -x_{(i_j)+1}h_{i_j} \]
and thus,
\[\begin{split}
c_{i_{j+1}} 
&= \{X + H_{i_{j+1}}, X + H_{i_{j+1}} + h_{i_{j+1}}\}\\
&= \{X +  H_{i_j} -x_{(i_j)+1}h_{i_j} , X + H_{i_j} -x_{(i_j)+1}h_{i_j}  + h_{i_{j+1}}\}.\\
\end{split}\]

The choice of $x_{(i_j)+1} = 0$ ensures $X + g \notin c_{i_{j+1}}$ unless $g=H_{i_j} + h_{i_{j+1}}$. So suppose
\begin{equation}\label{gehh}
g=H_{i_j} + h_{i_{j+1}}.
\end{equation}
Then setting $x_{(i_j)+1} = |Z|-1$ implies
\[
c_{i_{j+1}} = \{X +  H_{i_j} + h_{i_j} , X + H_{i_j} +h_{i_j}  + h_{i_{j+1}}\}.
\]
 Then $X + g \notin c_{i_{j+1}}$ unless 
 \begin{equation}\label{cij+1}
 g=H_{i_j} + h_{i_j}+ h_{i_{j+1}}.
 \end{equation}  
 But Equations~\ref{gehh} and \ref{cij+1} together imply $h_{i_j}=0$, a contradiction. Thus they both cannot be true, and one of the values 0 or $|Z|-1$ for $x_{(i_j)+1}$ ensures that $X + g \notin C_{i_{j+1}}$

We have shown that there are choices for $x_0, \ldots, x_{i_k} \in \Z_{\ge 0}$ such that for any choice of $x_{(i_k)+1}, \ldots, x_{d}$, $X+g \notin C_{i_k}$.  Since $h_i = 0$ for all $i > i_k$, we can set $x_i=0$ for all $i > i_k$ to ensure that $X + g \notin C_i$ for all $i > i_k$. Thus we have choices of $x_0, \ldots, x_{d}$ such that $X + g \notin C_{d-1}$, which proves $\chi$ is not $d$-polychromatic.
\end{proofof}

\section{Computer search for polychromatic linear colorings}\label{code}

By the fundamental theorem of finitely generated abelian groups, given $\ell$ and $d$, to test whether there is a $d$-polychromatic linear $Q_\ell$-coloring $\chi$ with $M$ colors, we need only consider linear colorings $\chi : \Z_{\geq 0}^{\ell+1} \to Z$ where $Z$ is of the form $Z= \bigoplus^{n}_{i=0} \Z/m_i$ and $M = m_0 m_1 \cdots m_n$. Suppose for $0 \le j \le \ell$,
\[\chi(f_j) = (c_{0j}, c_{1j}, \ldots, c_{nj}),\]
%\[\chi(f_j) = (c_{0j} \pmod{m_0}, c_{1j} \pmod{m_1}, \ldots, c_{nj}\pmod{m_n}),\]
where, as in the proof of Theorem~\ref{linnotopt}, we define $f_j \in \Z^{\ell+1}_{\ge 0}$ to be the vector whose $(j+1)$st entry is a 1, and all other entries are 0. Let $C_\chi$ be the matrix of coefficients 
\[C_\chi = \begin{pmatrix}
c_{00} &c_{01}& c_{02} & \ldots & c_{0\ell}\\
c_{10} &c_{11}& c_{12}& \ldots & c_{1\ell}\\ 
\vdots & \vdots & \vdots & \ddots & \vdots\\ 
c_{n0} &c_{n1}& c_{n2}& \ldots& c_{n\ell}
\end{pmatrix}.\]
Then for any counting vector $\vec{v} = (v_0, \ldots v_\ell)$,
\[\begin{split}\chi(v_0, \ldots, v_\ell) &=\left(\sum\limits_{j=0}^{\ell} c_{0j}v_j, \sum\limits_{j=0}^{\ell} c_{1j}v_j, \ldots, \sum\limits_{j=0}^{\ell} c_{nj}v_j  \right)\\
&= C_\chi\vec{v},
\end{split}\]
and a coloring $\chi$ is characterized by its matrix of coefficients $C_\chi$. Note that for $0 \le i \le n$, all arithmetic in the $i$th coordinate is done$\pmod{m_i}$. Further, we may assume that any entry $c_{ij}$ in $C_\chi$ is in $\Z/m_i$, and thus for a given factorization $M=m_0m_1 \cdots m_n$, there are a finite number ($M^\ell$) of possible linear colorings $\chi$ to consider.

The function \textit{isValidColoring} takes as parameters $\ell$, $d$, and $M$, and tests whether there is a $d$-polychromatic linear $Q_\ell$-coloring of the hypercube with $M$ colors. 
The pseudocode in Figure~\ref{pseud} describes the function, and we comment on it below.  Note all lists are indexed to start from 0. 
\begin{figure}
\begin{algorithmic}[1]
\Function{isValidColoring}{$\ell$, $d$, $M$}
\ForAll{factorizations $M=m_0m_1 \cdots m_n$ with $m_i>1$}
\ForAll{colorings $\chi:\Z_{\geq 0}^{\ell+1} \to  \bigoplus^{n}_{i=0} \Z/m_i$ with matrix $C_\chi$}
\State ChiIsPoly $\gets 1$
	\ForAll{counting vectors $\vec{x} = (x_0=0,x_1, \ldots, x_{d-1},x_d=0) \in [M]^{d+1}$}
		\State $S\gets \emptyset$
	
	   	\State
		\ForAll{star locations $\vec{s}$, $0=s_0<s_1 < \cdots <s_\ell<s_{\ell+1}=d+1$}

			\For{$i := 0$ to $\ell$}
				\State $y_i\gets \sum\limits_{k=s_{i}}^{s_{(i+1)} - 1} x_k$ 
			\EndFor
			\State
	\ForAll{$\vec{w} \in \Z^{\ell+1}_{\ge 0}$ where for $0 \le i \le \ell$, $0\le w_i < s_{i+1}-s_i$}
				\State $\vec{v} \gets \vec{y}+ \vec{w}$
				\State add $C_\chi\vec{v}$ to set S
			\EndFor
		\EndFor
		
		\State
		\If{$|S| \neq M$} %NOTE: $|S| = M$. 
			\State ChiIsPoly $\gets 0$
		\EndIf
	\EndFor
	\If {ChiIsPoly = 1}
	\State \Return TRUE
	\EndIf
\EndFor
\EndFor
\State
\State \Return FALSE
\EndFunction
\end{algorithmic}
\caption{Pseudocode for the function isValidColoring}\label{pseud}
\end{figure}

\textbf{Line 2} loops through all factorizations of $M=m_0m_1 \cdots m_n$.\\
\textbf{Line 3} loops through all colorings $\chi:\Z_{\geq 0}^{\ell+1} \to  \bigoplus^{n}_{i=0} \Z/m_i$ using the given factorization of $M$.\\
\textbf{Line 4:} ChiIsPoly is an indicator variable for whether $\chi$ is polychromatic, initially set to 1.\\
\textbf{Line 5} loops through every possible counting vector, $\vec{x}$, for a $Q_d$. Since $Z=\bigoplus^{n}_{i=0} \Z/m_i$, and for all $i$, $m_i$ divides $M$, we need only consider values of $x_i$ between 1 and $M$. Without loss of generality we set $x_0 = x_d = 0$.  Inside this loop we are concerned only with a $Q_d$ with this counting vector, which we refer to as ``the $Q_d$.''\\
\textbf{Line 6:} $S$ is the set of colors of the $Q_\ell$'s in the $Q_d$, initially set to $\emptyset$.\\
\textbf{Line 8} loops through possible star placements of $\ell$ stars to specify a $Q_\ell$ in the $Q_d$.  The star placements $s_1< \cdots < s_\ell$ are between 1 and $d$ (among the $d$ stars in the $Q_d$) and the values $s_0 = 0$ and $s_{\ell + 1} = d + 1$ represent two additional ``dummy'' stars which simplify the notation of subsequent computations. 
\\\textbf{Lines 9-11} create $\vec{y}$, the minimum number of 1's in the counting vector for any $Q_\ell$ with star locations given by $\vec{s}$, based on the counting vector for the $Q_d$.
\\\textbf{Lines 13-16} apply the coloring $\chi$ to every $Q_\ell$ with star locations given by $\vec{s}$, and add the resulting colors to the set $S$.  The quantity $w_i$ accounts for the possible number of 1's that may be added between the $i$th and $(i+1)$st stars in a $Q_\ell$ with stars in positions given by  $\vec{s}$. Thus the quantity $v_i=y_i + w_i$ represents the number of 1's between the $i$th and $(i+1)$st stars in a particular $Q_\ell$, and the vector $\vec{v} = \vec{y} + \vec{w}$ is the counting vector for a particular $Q_\ell$. The color $\chi(\vec{v}) = C_\chi\vec{v}$ is added to $S$. 
\\\textbf{Lines 19-21:} Now that all colors of all $Q_\ell$'s in the $Q_d$ have been added to $S$, if all $M$ colors are not in $S$, this implies the $Q_d$ does not contain a $Q_\ell$ of every color. Thus $\chi$ is not $d$-polychromatic and we set ChiIsPoly to 0. Then we return to Line 3 and test the next $Q_d$.\\
\textbf{Lines 23-25:} If ChiIsPoly=1, then under the current coloring $\chi$ every $Q_d$ contains a $Q_\ell$ of every color.  Thus $\chi$ is $d$-polychromatic and the program returns TRUE.\\
\textbf{Line 29:} If we reach this point in the program then every coloring has been tested and none are $d$-polychromatic.  Thus the program returns FALSE.

\subsection{Results of computer search}\label{dere}

\begin{proofof}{Proposition~\ref{small}}
Theorem~\ref{Chen} implies that $p_{lin}^2(4) \ge 6$, $p_{lin}^2(5) \ge 9$, and $p_{lin}^3(5) \ge 6$.  For the remaining three cases, in Figure~\ref{colorings} we exhibit a 6-polychromatic linear $Q_2$-coloring $\chi_{2,6}$ with 16 colors, a  6-polychromatic linear $Q_3$-coloring $\chi_{3,6}$ with 12 colors, and a  6-polychromatic linear $Q_4$-coloring $\chi_{4,6}$ with 9 colors.  (In the case $\ell=3$, $d=6$, there is more than one abelian group $Z$ where a coloring $\chi_{3,6}:\Z^{4}_{\ge 0} \to Z$ is 6-polychromatic, so we give multiple examples.) These colorings were found to be polychromatic using the function isValidColoring.

Theorem ~\ref{OS} implies $p^\ell(d) \le \binom{d+1}{\ell + 1}$.  Thus to prove these values are optimal, it suffices to examine all possible linear colorings onto all abelian groups with at most $\binom{d+1}{\ell + 1}$ elements.  We did this by evaluating the function isValidColoring for all $M \le \binom{d+1}{\ell + 1}$ for each $\ell$ and $d$.  For example, to establish $p_{lin}^2(4)=6$ we tested all linear $Q_2$-colorings with at most $\binom{5}{3} = 10$ colors and found no 4-polychromatic colorings with 7 or more colors. 
\begin{figure}
\begin{tabular}{l|l|l}
$\ell$  & Z & $\chi_{\ell,6}:\Z^{\ell+1}_{\ge 0} \to Z$ \\ \hline
2 & $\Z/4 \oplus \Z/4$ & $\chi_{2,6}(p,q,r) = (q-r \pmod{4}, p-r \pmod{4})$\\
3 & $\Z/3 \oplus \Z/4$ & $\chi_{3, 6}(p, q, r, s) = (r - s \pmod{3}, p-q+r-s \pmod{4})$\\
3 & $\Z/2 \oplus \Z/2 \oplus \Z/3$ & $\chi_{3, 6}(p, q, r, s) = (s \pmod{2}, q + r \pmod{2}, p+q-r \pmod{3})$\\
4 & $\Z/3 \oplus \Z/3$ & $\chi_{4,6}(p,q,r,s,t) = (q + s + t\pmod{3}, q-r+s\pmod{3})$\\
\end{tabular}
\caption{Optimal colorings for Proposition~\ref{small}}\label{colorings}
\end{figure}
\end{proofof}

\end{document}